\documentclass[12pt,twoside]{amsart}

\usepackage{amssymb,latexsym,bbm,graphicx,epsfig,epic,eepic,oldgerm,psfrag}
\usepackage{a4wide}

\usepackage{xypic} 
\input xy

\xyoption{all}

\theoremstyle{plain}
\newtheorem{theorem}{Theorem}[section]
\newtheorem{lemma}[theorem]{Lemma}
\newtheorem{proposition}[theorem]{Proposition}

\newtheorem*{remark*}{Remark}
\newtheorem*{remarks*}{Remarks}

\newtheorem{example}[theorem]{Example}

\newtheorem*{example*}{Example}
\newtheorem*{examples*}{Examples}

\newtheorem{systole.conjecture}[theorem]{Systole Conjecture}

\newtheorem*{notation*}{Notation}   
\newtheorem*{openproblem*}{Open Problem}

\newcommand{\proofend}{\hspace*{\fill} $\Box$\\}
\newcommand{\diam}{\hspace*{\fill} $\Diamond$}

\newcommand{\Tcirc}{\overset%
{\raisebox{-.3ex}[0ex][-.3ex]{\mbox{$\scriptscriptstyle\circ$}}\mskip-5mu}T}
\newcommand{\ecirc}{\overset%
{\raisebox{-.3ex}[0ex][-.3ex]{\mbox{$\scriptscriptstyle\circ$}}\mskip-5mu}e}

\def\1{\:\!}
\def\2{\;\!}
\def\s{\smallskip}
\def\m{\medskip}

\def\Diffc0{\operatorname{Diff^c_0}}
\def\Symp{\operatorname{Symp}}

\def\Sympc0{\operatorname{Symp^c_0}}

\def\width{\operatorname{width\2}}
\def\inradius{\operatorname{inradius\2}}
\def\reg{\operatorname{reg}}

\def\EHZ{\operatorname{EHZ}}

\def\ga{\alpha}
\def\gb{\beta}
\def\gg{\gamma}

\def\gve{\varepsilon}
\def\gf{\varphi}

\def\go{\omega}

\def\cb{{\mathcal B}}

\def\cf{{\mathcal F}}

\def\cp{{\mathcal P}}

\def\CC{\mathbbm{C}}

\def\RR{\mathbbm{R}}

\def\ZZ{\mathbbm{Z}}

\def\pp{\partial}

\def\ni{\noindent}

\def\m{\medskip}

\def\proof{\noindent {\it Proof. \;}}

\begin{document}

\title[]{Shortest closed billiard orbits on convex tables}

\author{Naeem Alkoumi}
\thanks{NA partially supported by 
the research fellowship 2013.0061
granted by the Federal Department of Home Affairs FDHA 
of the Swiss government}
\address{
    Naeem Alkoumi,
    Institut de Math\'ematiques,
Universit\'e de Neuch\^atel}
\email{naeem.alkoumi@unine.ch}

\author{Felix Schlenk}  
\thanks{FS partially supported by SNF grant 200020-144432/1.}
\address{Felix Schlenk,
Institut de Math\'ematiques,
Universit\'e de Neuch\^atel}
\email{schlenk@unine.ch}

\date{\today}

\begin{abstract}
Given a planar compact convex billiard table~$T$,
we give an algorithm to find the 
shortest generalised closed billiard orbits on~$T$.
(Generalised billiard orbits are usual billiard orbits if $T$
has smooth boundary.)
This algorithm is finite if $T$ is a polygon
and provides an approximation scheme in general.
As an illustration, we show that the shortest generalised closed
billiard orbit in a regular $n$-gon~$R_n$ is 2-bounce for~$n \ge 4$,
with length twice the width of~$R_n$.
As an application we obtain an algorithm computing the Ekeland--Hofer--Zehnder
capacity of the four-dimensional domain $T \times B^2$ in the standard symplectic vector space~$\RR^4$.
Our method is based on the work of Bezdek--Bezdek in~\cite{BezBez09}
and on the uniqueness of the Fagnano triangle in acute triangles. 
It works, more generally, for planar Minkowski billiards.
\end{abstract}

\maketitle

\section{Introduction and main results}

Mathematical billiards is a fascinating topic, with 
an abundance of problems and results. 
Almost every mathematical theory can be illustrated by and applied to 
a problem in mathematical billiards, see~\cite{Gut12, Kat05, KT91, Tab05, VGS92}
for excellent surveys. 
Here, we study the most elementary problem one can ask: 
Describe the set of shortest closed billiard orbits and their lengths
on a planar convex billiard table.

By a planar convex billiard table we mean a compact convex set~$T$ in~$\RR^2$
with non-empty interior~$\Tcirc$.
The boundary~$\pp T$ may be smooth or not, and $T$ may be strictly convex or not.
An outward support vector at~$q \in \pp T$ is a vector~$\nu$ such that
$$
\langle x - q, \nu \rangle \le  0 \quad \mbox{ for all } x \in T .
$$
A point~$q \in \pp T$ is called {\it smooth}\/ if the outward support vector of~$T$ at~$q_i$ in unique.
Equivalently, there is a unique line through~$q$ that is disjoint from~$\Tcirc$.

If $\pp T$ is smooth, a billiard orbit in~$T$ is a polygonal curve in~$T$
with vertices on~$\pp T$, such that at each vertex the incidence angle is equal to the reflection angle. 
Following~\cite{BezBez09, Gho04} we define a 
{\it generalised billiard orbit}\/ on~$T$ to be a sequence of points $q_i \in \pp T$, $i \in \ZZ$,
such that for every~$i$,
$$
\nu_i \,:=\, \frac{q_i - q_{i-1}}{\| q_i - q_{i-1} \|} + \frac{q_i - q_{i+1}}{\| q_i - q_{i+1} \|}
$$
is an outward support vector of~$T$ at~$q_i$.
We call the points $q_i$ the {\it bounce points}\/ of the generalised billiard orbit.
A billiard orbit is called {\it regular}\/ if all its bounce points are smooth, and singular otherwise.
If $\pp T$ is smooth, then the generalised billiard orbits on~$T$ are simply the billiard orbits on~$T$.

A generalised billiard orbit~$c$ is {\it closed}\/ or {\it periodic}\/
if there exists $n \ge 2$ such that $q_{i+n} = q_i$ for all~$i \in \ZZ$.
The smallest $n$ that works is the {\it period}\/ of~$c$, 
which is then called an $n$-bounce billiard orbit.
We throughout identify closed billiard orbits with the same trace.

\begin{example*}
{\rm
On a equilateral triangle, there are three 2-bounce orbits (that are singular), and
two 3-bounce orbits, the regular equilateral orbit and the singular orbit running along the boundary.

\begin{figure}[h]
 \begin{center}
  \leavevmode\epsfbox{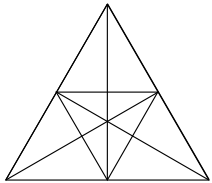}
 \end{center}
\caption{The five generalised closed billiard orbits on the equilateral triangle.}
 \label{figure.equilateral}
\end{figure}
%
%
}
\end{example*}

The length of an $n$-bounce orbit is of course defined by
$$
\ell (c) \,:=\, \sum_{i=0}^{n-1} \| q_{i+1}-q_i\| .
$$
\begin{notation*}
{\rm 
It will be convenient to use the following notation.
$$
\begin{array}{ll}
\cp (T) : & 
\mbox{ the generalised closed billiard orbits on~$T$} 
\\ [0.2em]
\cp_n (T) : & 
\mbox{ the $n$-bounce orbits in $\cp (T)$}
\\ [0.2em]
\cp_{\reg}(T) : &
\mbox{ the regular closed billiard orbits on~$T$}
\\ [0.2em]
\cp_{n,\reg} (T) : &
\mbox{ the $n$-bounce orbits in $\cp_{\reg} (T)$}
\\ [0.2em]
\cp_{\min} (T) : & 
\mbox{ the orbits in~$\cp (T)$ of minimal length.}
\end{array}
$$
}
\end{notation*}

Including singular orbits into the picture has many advantages. 
One advantage is the variational characterisation of~$\cp_{\min}(T)$ by Bezdek--Bezdek, 
that we recall in~Section~\ref{ss:variational}.
Another one is that generalised closed billiard orbits always exist.\footnote{While it is unknown 
whether every convex billiard table carries a regular closed orbit.
In fact, this is unknown even for general obtuse triangles.}
For instance there is a 2-bounce orbit of length $2 \2 \width (T)$,
where the width of~$T$ is the thickness of the thinnest band containing~$T$.
We can thus define
$$
\ell (T) \,:=\, \min \left\{ \ell (c) \mid c \in \cp (T) \right\} .
$$
The inradius of~$T$ is the radius of the largest disc contained in~$T$.
Ghomi proved in~\cite{Gho04} that always
\begin{equation} \label{e:42}
4 \inradius (T) \,\le\, \ell (T) \,\le\, 2 \2 \width (T) 
\end{equation}
with sharp lower bound if and only if $2 \inradius (T) = \width (T)$,
in which case $\cp_{\min}(T) \subset \cp_2(T)$.
Since $\width (T) \le 3 \inradius (T)$ for any convex set $T \subset \RR^2$ 
(see~\cite[Theorem~50]{Egg58}),
the bounds~\eqref{e:42} for $\ell (T)$ are sharp up to the factor~$\frac 32$.

In this note we describe a combinatorial process to find all shortest generalised
billiard orbits on a planar convex billiard table~$T$
and hence also~$\ell (T)$.
We outline the algorithm here. Details are given in Section~\ref{s:algo}.

\subsection{The algorithm}
Our algorithm is based on the following result of Bezdek--Bezdek from~\cite{BezBez09}:
\begin{equation} \label{e:inc23}
\cp_{\min}(T) \,\subset\, \cp_2(T) \cup \cp_{3,\reg}(T) .
\end{equation}
%
%
Assume first that $T$ is a polygon. 
The 2-bounce orbits, and in particular the 2-bounce orbits of minimal length $2 \width (T)$, 
are readily found.
In order to determine the shortest regular $3$-bounce orbits we recall that
in a triangle~$\Delta$ there is such an orbit if and only if $\Delta$ is acute, 
in which case this orbit is the Fagnano orbit, obtained by connecting the feet of the three altitudes of~$\Delta$.
If a polygon~$T$ has more than three edges, any regular $3$-bounce orbit on~$T$ is then found
as the Fagnano orbit of a triangle cut out by the lines supported by three edges of~$T$.
This leads to a finite algorithm for finding $\cp_{\min} (T)$ and $\ell (T)$,
that can be executed on a computer.

If $T$ is not polygonal, we approximate $T$ by a sequence of polygonal domains~$T_n$. 
Since $\ell (T)$ is continuous in the Hausdorff topology, $\ell (T_n)$ converges to~$\ell (T)$.
Moreover, if we take for each $n$ a shortest orbit $c_n \in \cp_{\min}(T_n)$, 
then a subsequence of~$c_n$ converges to an orbit $c \in \cp_{\min}(T_n)$,
and every orbit in $\cp_{\min} (T)$ can be obtained in this way.

\m
Several problems on closed orbits on planar convex billiard tables
are easier for tables with smooth boundary than for polygons. 
For instance, Birkhoff's famous theorem from~\cite{Bir27} asserts that
every strictly convex billiard table with smooth boundary carries 
infinitely many distinct closed orbits,
while for general polygons the existence of a regular closed orbit is unknown.
In contrast, our method uses a combinatorial process on polygons to give information on
shortest closed orbits on general convex billiard tables.

\subsection{Applications}\

\m \ni
{\bf 1. Some examples.}
To illustrate our method, we compute $\cp_{\min}(T)$ and $\ell (T)$
for triangles, for two classes of 4-gons, and for regular~$n$-gons, 
see Section~\ref{s:examples}.
For instance, the above algorithm immediately yields

\begin{proposition}
Let $R_n$ be a regular $n$-gon with $n \ge 5$ that is inscribed in the unit circle.
Then $\cp_{\min}(R_n) = \cp_2(R_n)$ and $\ell (R_n) = 2 \width (R_n) = 2 \left(1+ \cos \frac \pi n \right)$.
\end{proposition}

\ni
{\bf 2. 2-bounce orbits versus 3-bounce orbits.}
A problem posed by Zelditch asks whether the shortest billiard orbits on~$T$
are $2$-bounce or $3$-bounce (or both).
Our algorithm can decide this for polygons
and also for some non-polygonal convex billiard tables, see Section~\ref{s:zelditch}.

\m \ni
{\bf 3. Computation of the Ekeland--Hofer--Zehnder symplectic capacity.}
Endow $\RR^4$ with its standard symplectic form $\go_0 = dq_1 \wedge dp_1 + dq_2 \wedge dp_2$.
Denote by $B^4$ the open ball of radius~1 and by $Z^4$ the symplectic cylinder
$B^2 \times \CC$.
Let $\Symp (\RR^4)$ be the group of diffeomorphisms of~$\RR^4$ that preserve the 
symplectic form~$\go_0$.

\m 
A {\it symplectic capacity}\/ on~$(\RR^4, \go_0)$ associates with each subset~$S$ 
of~$\RR^4$ a number $c(S) \in [0,\infty]$ such that the following axioms are satisfied.

\m
{\bf (Monotonicity)}\;
$c (S) \le c (S')$ \, if\, $\gf (T) \subset T'$ for some $\gf \in \Symp (\RR^4)$;

\m
{\bf (Conformality)}\;
$c (r \1  T) = r^2  \2 c (T)$ \, for all\, $r>0$.

\m
{\bf (Nontriviality)}\;
$0 < c( B^4 )$ \, and \, $c( Z^4 ) < \infty$. 

\m
There are many different symplectic capacities, 
reflecting dynamical, geometric or holomorphic properties of a set
(see~\cite{CHLS07} for a survey).
The fascinating thing about capacities is that (in)equalities among them
imply relations between the different aspects of ``symplectic sets''.
Two dynamically defined symplectic capacities are the 
Ekeland--Hofer capacity and the Hofer--Zehnder capacity,
\cite{EkHo89, HZ90, HZ94}.
They agree on convex sets~$K$. 
Following~\cite{AO08} we denote their common value by $c_{\EHZ}(K)$.
Denote by $D^*T$ the unit ball bundle $T \times B^2 \subset \RR^2(q) \times \RR^2(p)$ in the cotangent bundle of~$T$.

\begin{proposition} \label{p:Ekeland}
For every compact convex set~$T \subset \RR^2$ it holds that
$c_{\EHZ}(D^*T) = \ell (T)$.
\end{proposition}

\proof
Monotonicity and conformality imply that $c_{\EHZ}$ is continuous in the Hausdorff topology. 
The same holds true for the function~$\ell$ in view of its
monotonicity and conformality property, see the end of Section~\ref{s:tools}.
We may thus assume that $T$ has smooth boundary.
For such billiard tables, the proposition is a ``folklore theorem'' known since the~1990th.
A precise treatment was given, however, only in~\cite{AO14}.
There, it is shown (in arbitrary dimensions) 
that $c_{\EHZ}(D^*T)$ is the minimum of $\ell (T)$ and the length of the shortest ``glide orbit''.
On a planar smooth convex billiard table~$T$, a glide orbit is simply an orbit running 
along the boundary~$\pp T$. Its length is thus larger than $2 \width (T)$.
Since $\ell (T) \le 2 \width (T)$ we conclude that $c_{\EHZ}(D^*T) = \ell (T)$.
\proofend

Symplectic capacities are very hard to compute in general.
In view of Proposition~\ref{p:Ekeland} our algorithm for computing 
$\ell (T)$ provides an algorithm for computing the capacity $c_{\EHZ}(D^*T)$,
finite if $T$ is polygonal and approximate otherwise.
For instance, for a regular $n$-gon with $n \ge 5$ odd we find 
$c_{\EHZ}(D^*R_n) = 2 \left(1+ \cos \frac \pi n \right)$.

\m
We conclude with addressing two problems.

\begin{itemize}
\item[1.] 
Is there an analogous algorithm for finding the shortest closed billiard orbits on tables of dimension $\ge 3$\1?

\s
\item[2.]
Does the algorithm also work for anisotropic billiards?
\end{itemize}

\s
{\bf Ad 1.}
Many works on billiards, such as~\cite{BezBez09, Gho04},
deal with convex domains of arbitrary dimension. 
Our method, however, seems to work only in dimension~two. 
Indeed, one of our main tools is the uniqueness of the Fagnano billiard orbit 
in acute triangles, and this result has no analogue in higher dimensions.

\s
{\bf Ad 2.}
While in this introduction we restricted ourselves to Euclidean billiards,
our method extends to anisotropic billiards, so-called Minkowski billiards.
In this generalisation of planar Euclidean billiards, 
there is given a (possibly non-symmetric) strictly convex body~$K \subset \RR^2$ with smooth boundary,
that determines the length of (oriented) straight segments and a reflection law for billiard orbits on~$T$.
The inclusion~\eqref{e:inc23} for generalised closed $K$-billiard orbits then still holds true,
but determining all shortest $2$-bounce and regular $3$-bounce orbits is somewhat harder, 
see Section~\ref{s:Minkowski}.
Since again
$c_{\EHZ} (T \times K)$ is the $K$-length of shortest generalised closed billiard orbits on~$T$,
we obtain an algorithm for computing the Ekeland--Hofer--Zehnder capacity of domains in~$\RR^4$
of the form $T \times K$ with $T,K \subset \RR^2$ convex.

\m
\ni
{\bf Acknowledgments.}
We wish to thank Yaron Ostrover and Sergei Tabachnikov for valuable discussions and suggestions.
We are particularly grateful to Sergei for showing us a proof of Lemma~\ref{le:FagMin}.
NA cordially thanks the Institut de Math\'ematiques for its hospitality in the academic year 2013--2014. 
He also thanks his family for the sacrifice it made to make this stay possible.
The present work is part of the author's activities within CAST, 
a Research Network Program of the European Science Foundation.

\section{Tools} \label{s:tools}

In this section we first recall two lemmata on 3-bounce billiard orbits in triangles, 
that we use to describe regular 3-bounce billiard orbits in convex polygons. 
We then rephrase the variational characterisation of shortest generalised closed billiard orbits
found by Bezdek--Bezdek.

\subsection{3-bounce billiard orbits}
 
A triangle is {\it acute}\/ if all its angles are $< \frac \pi 2$, 
it is {\it rectangular}\/ if one angle is~$\frac \pi 2$, and it is {\it obtuse}\/ if
one angle is $>\frac \pi 2$.

\begin{figure}[h]
 \begin{center}
  \leavevmode\epsfbox{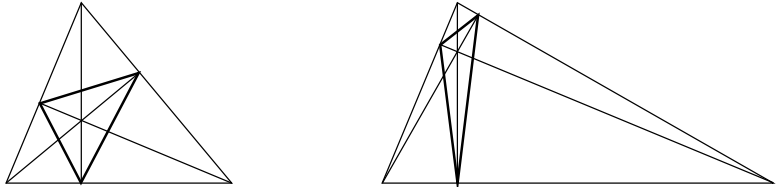}
 \end{center}
 \caption{Two Fagnano triangles}
 \label{figure.Fagnano}
\end{figure}

\s
Given an acute triangle~$T$,
the Fagnano triangle~$T_F$ of~$T$
is the triangle whose vertices are the feet of the three altitudes of~$T$,
see Figure~\ref{figure.Fagnano}.
It is named after J.\ F.\ de Tuschis a Fagnano, who around 1775 showed that
this triangle is the unique shortest triangle inscribed in~$T$,
and who also observed that this triangle represents a billiard orbit in~$T$.
For nice geometric proofs by Fej\'er and Schwarz see~\cite[\S 1.8]{Cox69} and~\cite[VII, \S 4]{CouRob79}.
These proofs, or a direct argument~\cite[p.\ 350]{CouRob79},
also show that the Fagnano triangle is shorter than twice the three altitudes of~$T$.
Another proof of uniqueness, that also applies to Minkowski billiards, is given in Lemma~\ref{le:FagMin}.

%
%

We begin with two well-known lemmata (see e.g.\ Proposition~9.4.1.3 in~\cite{Ber87}).

\begin{lemma} \label{le:acute}
Let $T$ be a triangle containing a regular $3$-bounce billiard orbit.
Then $T$ is acute.
\end{lemma}

\begin{figure}[h]
 \begin{center}
 \psfrag{u}{$u$}
 \psfrag{v}{$v$}
 \psfrag{w}{$w$}
 \psfrag{a}{$\ga$}
 \psfrag{b}{$\gb$}
 \psfrag{c}{$\gg$}
  \leavevmode\epsfbox{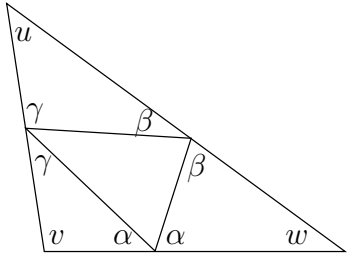}
 \end{center}
 \caption{The proof of Lemma~\ref{le:acute}}
 \label{figure.acute}
\end{figure}

\proof
Let $c$ be a regular $3$-bounce billiard orbit in~$T$, as in Figure~\ref{figure.acute}.
Then $\pi = \ga + \gb + \gg$, and 
$$
\pi \,=\, \ga + \gb + w \,=\, \gb  + \gg + u \,=\, \gg + \ga + v.
$$
Hence $u = \ga$, $v = \gb$, $w=\gg$, and therefore
$\pi > 2 \ga = 2u$, $\pi > 2v$, $\pi > 2w$, i.e., $T$ is acute.
\proofend

%
%

\begin{lemma} \label{le:unique}
Let $T$ be an acute triangle.
Then $T$ contains a unique regular $3$-bounce billiard orbit, forming the Fagnano triangle of~$T$.
\end{lemma}

\proof
Let $u,v,w$ be the angles of~$T$, and let $\Gamma$ be the triangle formed by a regular $3$-bounce billiard orbit
in~$T$.
As in the previous proof, $u = \ga$, $v = \gb$, $w=\gg$.
Hence $\ga' = \pi-2u$, $\gb' = \pi-2v$, $\gg' = \pi -2w$.
It follows that $T$ determines~$\Gamma$.
Hence $\Gamma$ is the Fagnano triangle.
\proofend

\begin{figure}[h]
 \begin{center}
 \psfrag{u}{$u$}
 \psfrag{v}{$v$}
 \psfrag{w}{$w$}
 \psfrag{a}{$\ga$}
 \psfrag{b}{$\gb$}
 \psfrag{c}{$\gg$}
 \psfrag{a'}{$\ga'$}
 \psfrag{b'}{$\gb'$}
 \psfrag{c'}{$\gg'$}
  \leavevmode\epsfbox{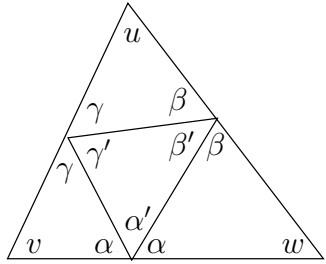}
 \end{center}
 \caption{The proof of Lemma~\ref{le:unique}}
 \label{figure.unique}
\end{figure}
%
%

\begin{proposition} \label{p:algo}
Let $T$ be a polygonal convex billiard table, and let $c$ be a regular $3$-bounce billiard orbit on~$T$.
Let $e_1, e_2, e_3$ be the edges of~$T$ hit by~$c$
(enumerated counterclockwise).
Denote by $\overline e_i$ the line supporting~$e_i$.
Then the lines $\overline e_1, \overline e_2, \overline e_3$ cut out an acute triangle~$\Delta$
containing~$T$, 
and the trace of~$c$ is the Fagnano triangle of~$\Delta$.
\end{proposition}

\proof
It is easy to see that $\overline e_1, \overline e_2$ are not parallel.
Since $T$ is convex, the point $\overline e_1 \cap \overline e_2$ lies on the right component of
$\overline e_1 \setminus \ecirc_1$.
Similarly, $\overline e_1 \cap \overline e_3$ lies on the left component of
$\overline e_1 \setminus \ecirc_1$.

\begin{figure}[h]
 \begin{center}
 \psfrag{v}{$v$}
 \psfrag{w}{$w$}
 \psfrag{a}{$\ga$}
 \psfrag{b}{$\gb$}
 \psfrag{c}{$\gg$}
 \psfrag{T1}{$e_1$}
 \psfrag{T2}{$e_2$}
 \psfrag{T3}{$e_3$}
 \psfrag{T1b}{$\overline e_1$}
 \psfrag{T2b}{$\overline e_2$}
 \psfrag{T3b}{$\overline e_3$}
  \leavevmode\epsfbox{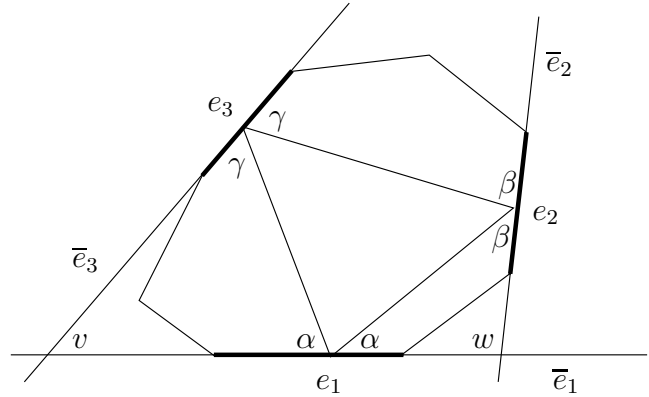}
 \end{center}
 \caption{The proof of Proposition~\ref{p:algo}}
 \label{figure.prop}
\end{figure}
%
%

\ni
With the angles as denoted in Figure~\ref{figure.prop} we have, as in the proof of Lemma~\ref{le:acute},
$$
v = \gb < \frac \pi 2, \quad w = \gg < \frac \pi 2 .
$$
Hence $v+w < \pi$. Hence $\overline e_1, \overline e_2, \overline e_3$ cut out a triangle~$\Delta$
enclosing~$T$ and with regular billiard orbit~$c$.
By Lemma~\ref{le:acute}, $\Delta$ is acute, and by Lemma~\ref{le:unique},
the trace of~$c$ is the Fagnano triangle of~$\Delta$.
\proofend

\subsection{The variational characterisation of shortest closed billiard orbits}
\label{ss:variational}

Let $T$ be a convex billiard table in~$\RR^2$.
Consider the set $\cb (T)$ of tuples $(q_1,q_2)$ and triples $(q_1,q_2,q_3)$ on the boundary~$\pp T$
that cannot be translated into the interior~$\Tcirc$.
Their length is defined as 
$$
\ell (q_1,q_2) = 2 |q_1-q_2|, \qquad \ell (q_1,q_2,q_3) = |q_1-q_2| + |q_2-q_3| +|q_3-q_1| .
$$
By compactness, $\ell = \min_{q \in \cb (T)} \ell (q)$ is attained.
Set 
$\cb_{\min} (T) = \left\{ q \in \cb (T) \mid \ell (q) = \ell \right\}$.

\begin{proposition} \label{p:BezBez}
{\rm (Bezdek--Bezdek~\cite{BezBez09})}
Let $T$ be a convex billiard table in~$\RR^2$.

\begin{itemize}
\item[(i)]
$\cp_{\min}(T) = \cb_{\min} (T)$.

\s
\item[(ii)]
A shortest generalised billiard orbit with 3 bounces is regular.
\end{itemize}
\end{proposition}

\proof
(i)
Let $\cf(T)$ be the set of 2-gons and 3-gons in~$\RR^2$ 
that cannot be translated into~$\Tcirc$.
Define two elements in~$\cf(T)$ to be equivalent if they are translates of each other.
It is shown in \cite[Lemma~2.4]{BezBez09} that the shortest elements of~$\cf(T)$,
up to equivalence, are the elements of~$\cp_{\min}(T)$.
Since the vertices of elements in~$\cp_{\min}(T)$ lie on~$\pp T$, 
each shortest equivalence class of~$\cf(T)$ contains an element of~$\cb_{\min} (T)$.

\s
(ii) 
If one of the vertices $q_1, q_2, q_3$~of $c$, say~$q_1$, is a non-smooth point of~$\pp T$, 
then it can be slightly moved along the boundary to a point $q_1'$ such that $(q_1',q_2,q_3)$ still
cannot be translated into the interior~$\Tcirc$ and so that the length of $(q_1',q_2, q_3)$ is less than $\ell (c)$; 
see the proof of Sublemma~3.1 in~\cite{BezBez09}.
\proofend

Denote again by~$\ell (T)$ the length of the orbits in $\cp_{\min} (T)$.
Proposition~\ref{p:BezBez}~(i) implies the following scale properties of~$\ell$.

\m
{\bf (Monotonicity)}\;
$\ell (T) \le \ell (T')$ \, if\, $T \subset T'$;

\s
{\bf (Conformality)}\;
$\ell (r \1 T) = r \2 \ell (T)$ \, for all\, $r>0$.

\m \ni
These two properties are useful for {\it estimating}\/ the shortest length~$\ell$:
If $\ell (S)$ is known and $S \subset T \subset r\1 S$, then monotonicity and conformality 
imply that
$$
\ell (S) \le \ell (T) \le r \2 \ell (S).
$$
For instance, assume that $S$ is a centrally symmetric convex billiard table with $S \subset T \subset r  S$.
By Corollary~1.3 in~\cite{Gho04} we have $\ell (S) = 2 \width (S)$.
Since the width is also monotone and conformal, we find that
$$
2 \width (S) \,\le\, \ell (T) \,\le\, 2 \width (T) \,\le\, 2 \2 r \width (S) .
$$
In the next section, we give an algorithm to {\it compute}\/ $\ell$.

\section{Algorithms} \label{s:algo}

Assume first that $T$ is a polygonal convex billiard table.
Proposition~\ref{p:algo} and~\ref{p:BezBez} give rise to finite algorithms for finding~$\cp_{\min}(T)$:
By Proposition~\ref{p:BezBez} we know that $\cp_{\min}(T) \subset \cp_2(T) \cup \cp_{3,\reg}(T)$. 
The set $\cp_2(T)$ is readily found, and $\cp_{3,\reg}(T)$ is found with the help of Proposition~\ref{p:algo}.
The set $\cp_{\min}(T)$ is then obtained by selecting the orbits of shortest length~$\ell (T)$.


\m \ni
{\bf Algorithm 1 (finding $\cp_2(T)$)}

\s \ni
{\it
\begin{itemize}
\item[(i)]
If $\pp T$ has parallel edges $e_i,e_j$, then the segments orthogonal to $\ecirc_i, \ecirc_j$
form regular $2$-bounce orbits on~$T$,
and all regular $2$-bounce orbits on~$T$ are of this form.

\s
\item[(ii)]
Given a vertex~$v$ and a disjoint edge~$e$, form the altitude~$s$ from $v$ to~$\overline e$.
Then $s$ is half of a generalised 2-bounce orbit on~$T$ if and only if
the end point of~$s$ lies on~$e$ and the line through~$v$ orthogonal to~$s$ is disjoint from~$\Tcirc$. 

\s
\item[(iii)]
Given two different vertices~$v_i, v_j$, 
the segment $s = v_iv_j$ is half of a generalised 2-bounce orbit on~$T$ if and only if
the lines through~$v_i,v_j$ orthogonal to~$s$ are disjoint from~$\Tcirc$.
\end{itemize}
}

\s
If one is only interested in finding the shortest $2$-bounce orbits, namely those of length $2 \width (T)$,
it suffices to look at the orbits arising in points (i) and (ii), since those in (iii) that are 
not covered by (ii) are longer.
Similarly, if one is only interested in finding~$\ell (T)$, it suffices to look at the orbits arising in (ii).
By Proposition~\ref{p:algo} we have

\m \ni
{\bf Algorithm 2 (finding $\cp_{3,\reg}(T)$)}

\s \ni
{\it
Take all triples $e_1, e_2, e_3$ among the edges of~$T$
that cut out an acute triangle containing~$T$.
Among these triangles, select those whose Fagnano triangle is contained in~$T$,
i.e., $\overline e_i \cap \overline e_j$ projects to~$\ecirc_k$ for $\{i,j,k\} = \{1,2,3\}$.
}

\m
As the two algorithms show, 
the lengths of the orbits in $\cp_2(T) \cup \cp_{3,\reg}(T)$ can be computed in terms of
the coordinates of the vertices of~$T$.
The whole algorithm can thus be executed by a computer code.

\m
We now use the above algorithms to investigate $\cp_{\min}(T)$ and~$\ell (T)$
on arbitrary planar convex billiard tables~$T$.
Let $T$ be such a table.
Fix $\gve >0$.
Choose a polygonal convex billiard table~$T_\gve$ such that
\begin{equation} \label{e:inc}
T_\gve \,\subset\, T \,\subset\, (1+\gve) \, T_\gve .
\end{equation} 
By monotonicity and conformality of~$\ell$,
\begin{equation} \label{e:sandwich}
\ell (T_\gve) \,\le\, \ell (T) \,\le\, (1+\gve) \, \ell (T_\gve) .
\end{equation}
Take a sequence $\gve_n \to 0$ and corresponding polygonal convex billiard tables~$T_{\gve_n}$
satisfying~\eqref{e:inc}.
For each~$n$ choose $c_n \in \cp_{\min}(T_{\gve_n})$.
Since each~$c_n$ has 2 or 3 bounces, Proposition~\ref{p:BezBez}~(i) and~\eqref{e:inc} imply
that a subsequence of~$c_n$ converges to an orbit $c \in \cp_{\min}(T)$.
On the other hand, it is clear that every $c \in \cp_{\min}(T)$ can be obtained in this way.

Summarizing, we see that given $\gve >0$ we have a finite algorithm computing a number $\ell_\gve (T)$
such that 
$$
\ell (T) \,\le\, \ell_\gve (T) \,\le\, (1+\gve) \, \ell (T) .
$$

\section{Examples} \label{s:examples}

To illustrate our method, we compute the shortest generalised closed billiard orbits in triangles,
in two special classes of 4-gons and in regular $n$-gons.
Throughout we apply Algorithm~2.

\subsection{Shortest billiard orbits in triangles}

\begin{proposition} \label{p:triangle}
Let $T$ be a triangle.

\begin{itemize}
\item[(i)]
If $T$ is acute, the shortest generalised closed billiard orbit on~$T$ is the Fagnano triangle
(which is regular).

\s
\item[(ii)]
If $T$ is rectangular or obtuse, the shortest generalised closed billiard orbit on~$T$ is 
the singular $2$-bounce orbit starting at the vertex with angle $\ge \frac \pi 2$.
In particular, $\ell (T) = 2 \width (T)$.
\end{itemize}
\end{proposition}

\proof
(i)
The shortest $2$-bounce orbits lie on (one or two or three of) the altitudes of~$T$,
and by Lemma~\ref{le:unique}, the Fagnano triangle is the only regular 3-bounce orbit.
It is shorter than twice the three altitudes of~$T$.

\m
(ii)
Let $h_1$ be the altitude starting at the vertex~$v_1$ with angle $\ge \frac \pi 2$.
If $T$ is rectangular, $h_1$ is shorter than the other two altitudes, which lie on the edges
containing~$v_1$.
If $T$ is obtuse, $h_1$ is the only altitude contained in~$T$.
By Lemma~\ref{le:acute}, $T$ contains no regular 3-bounce orbit.
\proofend

\subsection{Shortest billiard orbits in two special classes of 4-gons}

In this paragraph we look at convex 4-gons that either have two parallel edges or a rectangular corner.

\subsubsection{4-gons with two parallel edges}
\begin{figure}[h]
 \begin{center}
 \psfrag{a1}{$\ga_1$}
 \psfrag{a2}{$\ga_2$}
 \psfrag{e1}{$e_1$}
 \psfrag{e2}{$e_2$}
 \psfrag{e3}{$e_3$}
  \leavevmode\epsfbox{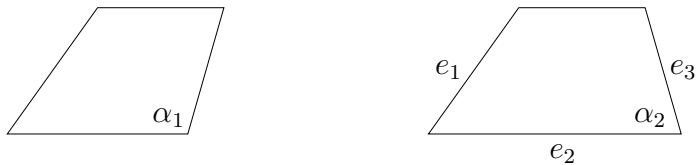}
 \end{center}
 \caption{Two 4-gons with two parallel edges}
 \label{figure.4-gon1}
\end{figure}

Up to isometry, such a polygon looks like one of the polygons in Figure~\ref{figure.4-gon1},
where $\ga_1 \ge \frac \pi 2$ and $\ga_2 < \frac \pi 2$.

%
%

\begin{figure}[h]
 \begin{center}
  \leavevmode\epsfbox{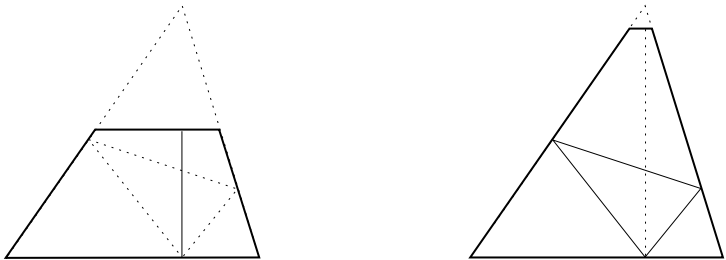}
 \end{center}
 \caption{Two possibilities for the orbits attaining $\ell (T)$}
 \label{figure.4-gon2}
\end{figure}

\ni
In the first case, there is no triple $e_1, e_2, e_3$ among the edges of~$T$
that cuts out an acute triangle containing~$T$.
Hence $\cp_{\min}(T) \subset \cp_2(T)$ and $\ell (T) = 2 \width (T)$.
In the second case, the only triple $e_1, e_2, e_3$ 
that cuts out an acute triangle containing~$T$ is as marked in Figure~\ref{figure.4-gon1}.
The Fagnano triangle~$\Delta_F$ of the corresponding triangle may lie in~$T$ or not.
If it does, then $\ell (T) = \min \left\{ \ell (\Delta_F), 2 \1 h \right\}$,
where~$h$ is the distance between the two parallel edges of~$T$.
Both possibilities for the minimum occur as Figure~\ref{figure.4-gon2} illustrates.

%
%

\subsubsection{4-gons with a rectangular corner}
In view of the previous example, we can assume that no two edges of~$T$ are parallel.
Since the angle sum is $2 \pi$,
$T$ then looks up to isometry like one of the following three polygons: 

%
\begin{figure}[h]
 \begin{center}
 \psfrag{a}{$a$}
 \psfrag{b}{$b$}
 \psfrag{c}{$c$}
 \psfrag{d}{$d$}
 \psfrag{aa}{$\ga$}
 \psfrag{bb}{$\gb$}
 \psfrag{cc}{$\gg$}
 \psfrag{1}{\mbox{(1): $\ga, \gg < \frac \pi 2$, $\gb > \frac \pi 2$}}
 \psfrag{2}{\mbox{(2): $\ga, \gg > \frac \pi 2$, $\gb < \frac \pi 2$}}
 \psfrag{3}{\mbox{(3): $\ga < \frac \pi 2$, $\gg > \frac \pi 2$}}
 \psfrag{v}{$v$}
  \leavevmode\epsfbox{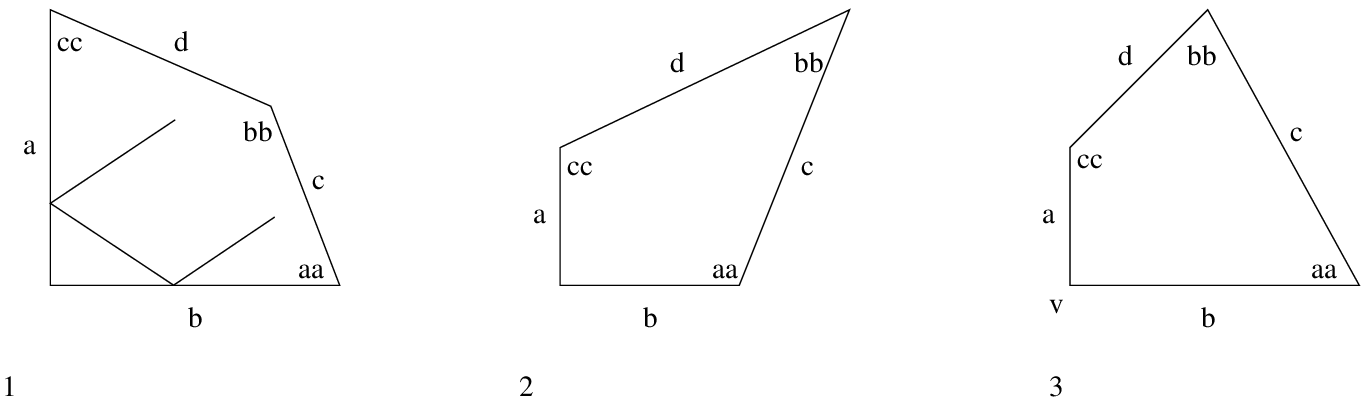}
 \end{center}
 \label{figure.ortho}
\end{figure}
%
%

\ni
In Case 3, $\gb$ may be acute or not.

\s \ni
There is no closed 3-bounce billiard orbit with bounces on~$a$ and~$b$
since for such an orbit two legs would be parallel (see the left picture of Figure~\ref{figure.ortho}).
A closed 3-bounce billiard orbit must thus bounce on $acd$ or $bcd$ (up to orientation).
In Cases~1 and~2, the triples $\overline a, \overline c, \overline d$ and 
$\overline b, \overline c, \overline d$ do not cut out an acute triangle containing~$T$,
and the same holds true in Case~3 for the triple $\overline a, \overline c, \overline d$
and if $\beta \ge \frac \pi 2$ also for the triple $\overline b, \overline c, \overline d$.
Hence $\cp_{\min}(T) \subset \cp_2(T)$ and $\ell (T) = 2 \width (T)$.
In Case~3 with $\beta$ acute, the Fagnano triangle~$\Delta_F$ of~$\overline b, \overline c, \overline d$
may or may not lie in~$T$. 
If it does, then $\ell (T) = \min \left\{ \ell (\Delta_F), 2 \1 h \right\}$,
where~$h$ is the distance from $v$ to~$c$.
Again, both possibilities for the minimum occur.

\subsection{Shortest billiard orbits in regular $n$-gons}

For $n \ge 3$ consider a regular $n$-gon $R_n$.
For $n$ even, $R_n$ is centrally symmetric, 
and hence $\cp_{\min} (R_n) \subset \cp_2(R_n)$ by Corollary~1.3 in~\cite{Gho04}.
This holds true for all $n \ge 4$.
More precisely, we have

\begin{proposition} \label{p:Tn}
Consider the regular $n$-gon~$R_n$ inscribed in the circle of radius~$1$.

\s
\begin{itemize} 
\item[(i)]
If $n=3$, then $\cp_{\min} (R_n)$ consists of the Fagnano orbit of~$T_3$.
Its length is $\frac{3 \sqrt 3}{2}$.

\m
\item[(ii)]
If $n \ge 5$ is odd, then $\cp_{\min} (R_n)$ consists of the $n$ singular 2-bounce orbits
starting at the vertices of $R_n$. Their length is $2 \left( 1+\cos \frac{\pi}{n} \right) = 2 \width (R_n)$. 

\m
\item[(iii)]
If $n$ is even, then $\cp_{\min} (R_n)$ consists of the $\frac n2$ bands of 2-bounce orbits of length 
$4 \cos \frac \pi n =2 \width (R_n)$.

\m
\item[(iv)]
If $n = 3k$, then there exist $k$ regular $3$-bounce orbits on~$R_n$, namely the equilateral triangles with vertices
on the midpoints of the edges they hit. 
Their length is $3 \sqrt{3} \cos \frac \pi n$ which is larger than $2 \width (R_n)$ if $k \ge 2$.
If $n \neq 3k$, then $\cp_{3,\reg}(R_n)$ is empty.
\end{itemize}
\end{proposition}

\proof
The length of an edge of~$R_n$ is $\bigl| 1-e^{\frac{2\pi i}{n}} \bigr| = 2 \sin \frac{\pi}{n}$.
Hence the distance between the origin and the midpoint of an edge is $\cos \frac \pi n$,
and so
$$
\width (R_n) \,=\, 
\left\{\begin{array}{ll}
2 \cos \frac{\pi}{n}  &  \mbox{if }\;  \mbox{$n$ is even}, \\ [.3em]
1+\cos \frac{\pi}{n}  &  \mbox{if }\;  \mbox{$n$ is odd}. 
\end{array}\right.
$$

(i) is Proposition~\ref{p:triangle}.
The 2-bounce orbits on $R_n$ are obvious. 
It remains to determine all regular $3$-bounce orbits on~$R_n$ for $n \ge 4$.

Let $c \in \cp_{3,\reg}(R_n)$, with bounce points on the edges $e_{i_1}, e_{i_2}, e_{i_3}$.
Assume first that $n=3k$ and that $\{ i_1, i_2, i_3 \}$ is of the form $\{ i, i+k, i+2k \}$.
Then $\overline e_i, \overline e_{i+k}, \overline e_{i+2k}$ cut out an equilateral triangle~$\Delta$
containing~$R_n$. 
By Lemma~\ref{le:unique}, $c$ runs along the Fagnano triangle of~~$\Delta$, which is equilateral.

\begin{figure}[h]
 \begin{center}
 \psfrag{0}{$0$}
 \psfrag{1}{$1$}
 \psfrag{L}{$L$}
 \psfrag{d}{$d$}
 \psfrag{e1}{$e_{i_1}$}
 \psfrag{e2}{$e_{i_2}$}
 \psfrag{e3}{$e_{i_3}$}
  \leavevmode\epsfbox{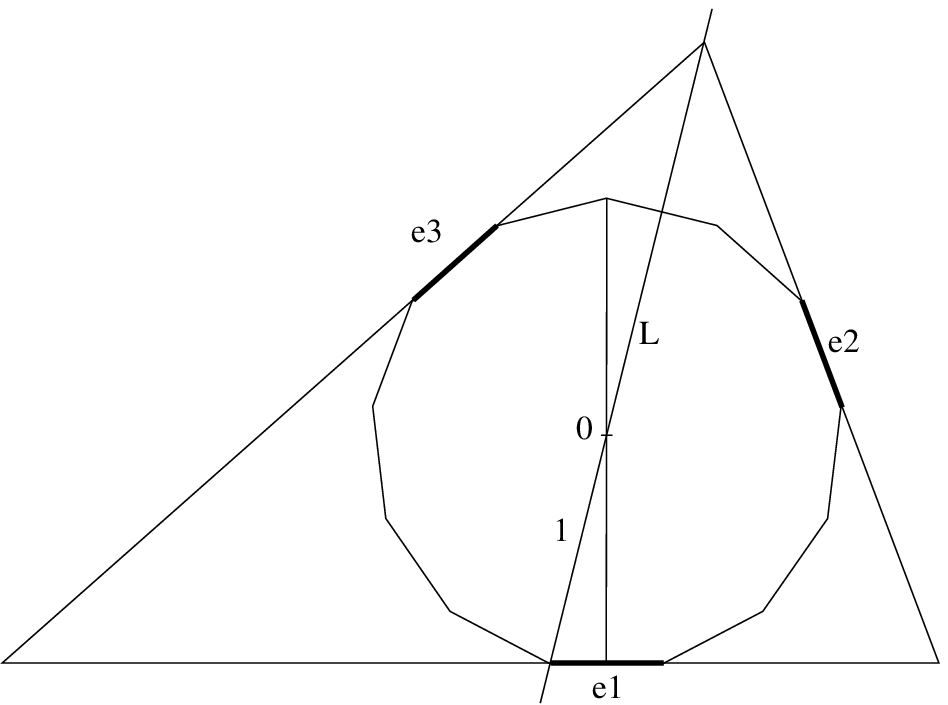}
 \end{center}
 \label{figure.Tn}
\end{figure}
%
%

Assume now that $n \neq 3k$ or that $n=3k$ and $\{ i_1, i_2, i_3 \}$ is not of the form $\{ i, i+k, i+2k \}$.
Since $c \in \cp_{3,\reg}(R_n)$, the lines $\overline e_{i_1}, \overline e_{i_2}, \overline e_{i_3}$ cut out an equilateral triangle~$\Delta$ containing~$R_n$. 
We must show that the Fagnano triangle~$\Delta_F$ of~$\Delta$ is not contained in~$R_n$.
Assume first that $n$ is odd.
Denote by $\rho_y$ the reflection along the $y$-axis.
After renaming $i_1, i_2, i_3$, if necessary, we can assume that $e_{i_1}, e_{i_2}, e_{i_3}$
are as in 
the figure:
$e_{i_1}$ is the lower horizontal edge, and $\rho_y (e_{i_2}) \neq e_{i_3}$, with $e_{i_2}$ below $e_{i_3}$. 
The hardest case is when $e_{i_3}$ neighbors $\rho_y (e_{i_2})$, as in the figure. 
Then the line~$L$ through the vertex $v$ of~$\Delta$ and through~$0$ passes through the left boundary point of~$e_{i_1}$.
Hence a point~$q$ on~$L$ projects to~$e_{i_1}$ if and only if $|q| \le 1$.
Since $|v| >1$, we see that~$v$ does not project to~$e_{i_1}$. 
Hence $\Delta_F$ is not contained in~$R_n$.
If $e_{i_3}$ does not neighbor $\rho_y (e_{i_2})$, then $v$ will project to a point on~$\overline e_{i_1}$
even further apart from~$e_{i_1}$.
The argument for $n$ even is similar and left to the interested reader.
\proofend

\section{Application to a question of Zelditch}
\label{s:zelditch}

Let again $T$ be a planar convex billiard table, and recall from Proposition~\ref{p:BezBez} that 
$$
\cp_{\min}(T) \,\subset\, \cp_2(T) \cup \cp_{3,\reg}(T) .
$$
It is interesting to ask when $\cp_{\min}(T) \subset \cp_2(T)$.
This problem was brought up by Zelditch~\cite{Zel00} in relation with  
the inverse spectral problem on smooth domains.

For polygonal convex billiard tables, our algorithm solves this problem, 
cf.\ the examples in the previous section.
Classes of convex billiard tables with $\cp_{\min}(T) \subset \cp_2(T)$ are 
centrally symmetric tables or, more generally, tables with $2 \inradius (T) = \width (T)$,
see ~\cite{Gho04}, 
and so-called fat disc-polygons~\cite{BezBez09}.

Non-polygonal examples with $\cp_{\min}(T) \subset \cp_{3,\reg}(T)$ can be obtained as follows:
Let $T$ be a convex billiard table 
and assume that there exists $c \in \cp_3(T)$ with $\ell (c) < 2 \width (T)$.
Then for any convex billiard table~$T'$ with
$$
r_1 \, T \subset T' \subset r_2 \, T \qquad \mbox{ and } \qquad
\frac{r_2}{r_1} \,<\, \frac{2 \width (T)}{\ell (c)} 
$$
we still have $\cp_{\min}(T) \subset \cp_{3,\reg}(T)$.
Indeed, using monotonicity and conformality of~$\ell$ and of the width we can estimate
$$
\ell (T') \,\le\, r_2 \, \ell (T) \,\le\, r_2 \, \ell (c) \,<\,
2 \2 r_1 \, \width (T) \,\le\, 2 \2 \width (T') .
$$  
Since the shortest generalised 2-bounce orbits on~$T'$ have length $2 \1 \width (T')$,
the claim follows.

\begin{example}
{\rm 
For the equilateral triangle $\Delta$ of edge length~$1$, 
the Fagnano triangle is also equilateral, and has length $\frac 32 < 2 \width (\Delta) = \sqrt 3$.
Hence for any convex billiard table $T'$ with 
$$
r_1 \, \Delta \subset T' \subset r_2 \Delta \qquad \mbox{ and } \qquad
\frac{r_2}{r_1} < \frac{2 \sqrt 3}{3} 
$$
every shortest generalised billiard orbit is a regular 3-bounce orbit.
}
\diam
\end{example}

\section{Generalisation to planar Minkowski billiards}  \label{s:Minkowski}
Many newer works on (shortest) billiard orbits on convex domains $T \subset \RR^2$
treat the more general case of Minkowski billiards:
There is given a strictly convex body~$K \subset \RR^2$ with smooth boundary,
which is used to define the length of straight segments in~$\RR^2$ and a reflection law on~$T$,
see \cite{ABKS, AKO, AO14, GutTab02}.

\begin{figure}[h]
 \begin{center}
 \psfrag{a}{$a$}
 \psfrag{b}{$b$}
 \psfrag{x}{$x$}
 \psfrag{K}{$K$}
 \psfrag{L1}{$L_1$}
 \psfrag{L2}{$L_2$}
 \psfrag{pT}{$\pp T$}
 \psfrag{TpT}{$T_x \pp T$}
  \leavevmode\epsfbox{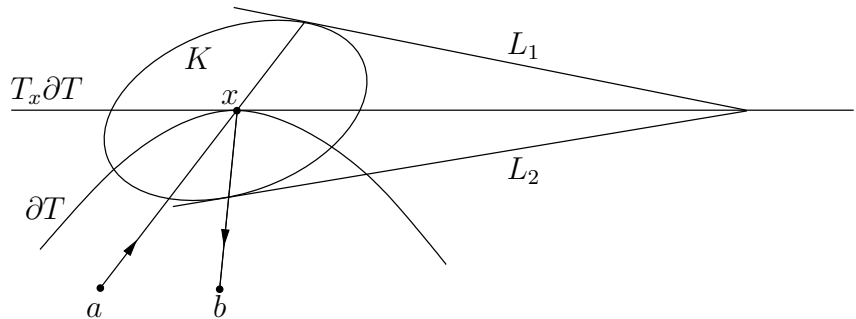}
 \end{center}
 \caption{The reflection law, geometrically}
 \label{figure.law}
\end{figure}

For symmetric $K$, the reflection law can be formulated as follows, \cite[\S 3]{GutTab02}.
Given interior points $a,b \in \Tcirc$ and a smooth boundary point $x \in \pp T$,
the segments $ax, xb$ are part of a $K$-billiard orbit on~$T$ if and only if $x$ is a critical point of the function
$y \mapsto \ell_K(ay) + \ell_K(yb)$ on~$\pp T$.
Equivalently, the exit direction~$xb$ can be found from the entrance direction~$ax$ and from~$K$ by drawing 
first the tangent line~$L_1$ and then the tangent line~$L_2$ to~$K$ as in Figure~\ref{figure.law}.

%
%

\ni
If $x$ is not smooth, then we agree that the reflection law holds at~$x$ if it holds with respect to some 
line that passes through~$x$ and is disjoint from~$\Tcirc$. 
For $K$ the unit disc, 
this reflection law and the associated billiard dynamics becomes the Euclidean one defined in the introduction. 
For the definition of the reflection law for non-symmetric~$K$ we refer to~\cite{ABKS, AKO, AO14}.
Note that for symmetric $K$, the length of a closed orbit does not depend on its orientation, 
but for non-symmetric~$K$ it may.

Our method extends to this more general setting. 
Indeed, as noticed in~\cite[\S 2.1]{ABKS},
the variational characterisation of~$\cp_{\min}(T,K)$ in Proposition~\ref{p:BezBez} 
still holds true in this setting.
In particular, 
the shortest generalised closed $K$-billiard orbits on~$T$ are 2-bounce or 3-bounce, 
and shortest 3-bounce orbits are regular.
It remains to find an efficient way to determine these orbits.
This is less straightforward than in the Euclidean case.

%

From now on we assume that $K$ is symmetric.
We first determine the set $\cp_2(T;K)$ of generalised 2-bounce $K$-billiard orbits on~$T$.
We start with a few observations.

\begin{itemize}
\item[(i)]
Given a regular 2-bounce orbit between edges $e_i,e_j$, these edges must be parallel by the symmetry of~$K$.
Moreover, by the strict convexity of~$K$, there is a unique band of parallel 2-bounce orbits 
between $\overline e_i, \overline e_j$. 

\s
\item[(ii)]
Given a point $v$ disjoint from a line~$L$, there is a unique point $v_L$ on~$L$ at which the $K$-distance
from $v$ to~$L$ is attained, because $K$ is strictly convex. 
We call the segment $vv_L$ the $K$-altitude from~$v$ to~$L$.

\s
\item[(iii)]
Given a segment $s$ there are unique parallels $L_1, L_2$
through the end points of~$s$ such that $s$ is a $K$-altitude from~$L_1$ to~$L_2$,
again because $K$ is strictly convex.
\end{itemize}

%
%
%

\s
With these observations, we obtain as in Section~\ref{s:algo} the following

\m \ni
{\bf Algorithm 1 (finding $\cp_2(T;K)$)}

\s \ni
{\it
\begin{itemize}
\item[(i)]
If $\pp T$ has parallel edges $e_i,e_j$, then the altitudes between $\overline e_i, \overline e_j$ 
that are based on $\ecirc_i, \ecirc_j$ form regular 2-bounce orbits on~$T$,
and all regular $2$-bounce orbits on~$T$ are of this form.

\s
\item[(ii)]
Given a vertex~$v$ and a disjoint edge~$e$, form the $K$-altitude~$s$ from $v$ to~$\overline e$.
Then $s$ is half of a generalised 2-bounce orbit on~$T$ if and only if
the end point of~$s$ lies on~$e$ and the line through $v$ parallel to~$\overline e$ is disjoint from~~$\Tcirc$. 

\s
\item[(iii)]
Given two different vertices~$v_i, v_j$, 
the segment $s = v_iv_j$ is half of a generalised 2-bounce orbit on~$T$ if and only if
the parallel lines $L_i, L_j$ through $v_i,v_j$ for which $s$ is a $K$-altitude
are disjoint from~$\Tcirc$.
\end{itemize}
}

\m
It remains to understand the regular 3-bounce orbits in Minkowski triangles~$\Delta$.
In~\cite{GutTab02} such triangles are called Fagnano triangles.

\begin{lemma} \label{le:FagMin}
Let $\Delta$ be a triangle in the Minkowski plane $(\RR^2, K)$.
Then there exists at most one Fagnano triangle in~$\Delta$.
\end{lemma}

\proof
The following proof was shown to us by Sergei Tabachnikov.
Given an oriented line~$L$ in~$\RR^2$ denote by $\sphericalangle L$ the oriented angle 
from the positively oriented $x$-axis to~$L$.
For $i=1,2$ let $u_i$ be an incoming billiard leg reflecting on a given line to~$v_i$.
Assume that $\sphericalangle u_1 > \sphericalangle u_2$, as in the left figure. 
Then the strict convexity of~$K$ implies that $\sphericalangle v_1 < \sphericalangle v_2$,
cf.\ Figure~\ref{figure.law}. 

\begin{figure}[h]
 \begin{center}
 \psfrag{u1}{$u_1$}
 \psfrag{u2}{$u_2$}
 \psfrag{v1}{$v_1$}
 \psfrag{v2}{$v_2$}
 \psfrag{P1}{$P_1$}
 \psfrag{P2}{$P_2$}
 \psfrag{Q1}{$Q_1$}
 \psfrag{Q2}{$Q_2$}
 \psfrag{R1}{$R_1$}
 \psfrag{R2}{$R_2$}
  \leavevmode\epsfbox{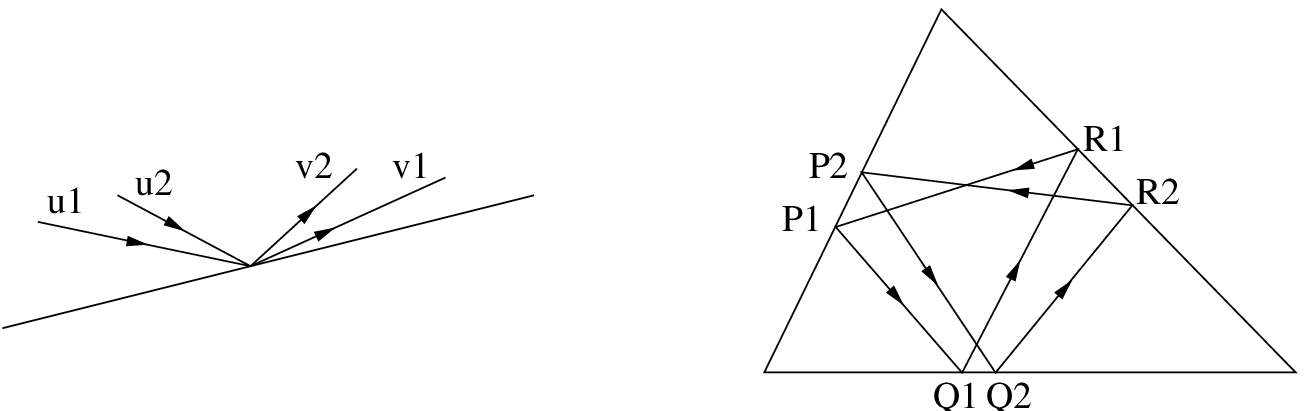}
 \end{center}
 \label{figure.Minmin}
\end{figure}
%
%

Now suppose that $P_1Q_1R_1$ and $P_2Q_2R_2$ are two different Fagnano triangles in~$\Delta$.
Then not all pairs of the respective sides of these triangles are parallel, 
say $\sphericalangle P_1Q_1 > \sphericalangle P_2Q_2$. 
Then $\sphericalangle Q_1R_1 < \sphericalangle Q_2R_2$, hence $\sphericalangle R_1P_1 > \sphericalangle R_2P_2$, 
hence $\sphericalangle P_1Q_1 < \sphericalangle P_2Q_2$, a contradiction. 
\proofend

The same argument shows that embedded $n$-bounce orbits in convex Minkowski $n$-gons are unique (if they exist).
Following~\cite{GutTab02} we call a triangle {\it acute}\/ if it admits a Fagnano orbit.
As in the Euclidean case we have

\m \ni
{\bf Algorithm 2 (finding $\cp_{3,\reg}(T,K)$)}

\s \ni
{\it
Take all triples $e_1, e_2, e_3$ among the edges of~$T$
that cut out an acute triangle containing~$T$.
Among these triangles, select those whose Fagnano triangle is contained in~$T$.
}

\m
Solving the following problem would complete the algorithm finding $\cp_{\min} (T;K)$ for symmetric~$K$. 

\begin{openproblem*}
Give an algorithm finding the Fagnano triangle in a Minkowski triangle.
\end{openproblem*}


\end{document}